\documentclass[11pt,letterpaper]{article}
\usepackage[hmargin=1in, vmargin=1in]{geometry}
\usepackage{graphicx}
\usepackage{amsmath}
\usepackage{dsfont}
\usepackage{indentfirst}
\usepackage{sidecap}
\usepackage{float}
\usepackage{multicol}
\usepackage{caption}
\usepackage{subcaption}
\usepackage[hyphens]{url}
\usepackage{hyperref}

\usepackage{fancyhdr}
\title{3D Printing A Pendant with A Logo}
\author{Prof. Edward Aboufadel \\ Department of Mathematics \\ Grand Valley State University}
\date{Version 1.0, July 2015}

\begin{document}
\maketitle

\section{Introduction}

The purpose of this short paper is to describe a project to manufacture a 3D-print of a pendant that includes a logo.  The methods described in this paper involve processing the image of the logo through a Mathematica script.  These methods can be applied to many logos and other images.  With the Mathematica script, a STereoLithography (.stl) file is created that can be used by a 3D printer.  Finally, the object is created on a 3D printer.   We assume that the reader is familiar with the basics of 3D printing.

\section{The Fighting Pancreases of Zachary University}

JDRF (formerly known as the Juvenile Diabetes Research Foundation) is a major charity devoted to preventing, treating, and curing type-1 diabetes.  Type-1 diabetes (also called T1D) occurs when the immune system is triggered to start destroying insulin-producing beta cells in the pancreas.  The cause of this trigger is uncertain, and a cure is not in sight at this time.  People with T1D introduce artificial insulin into their bodies through injections or an insulin pump, and continuous glucose monitors are now available in order to track glucose levels in the blood at five-minute intervals.  Type-1 diabetes should not be confused with the more common type-2 diabetes that afflicts a considerable number of Americans. About 1.25 million Americans have T1D.

There is a member of the author's family afflicted with T1D.  Each year in late summer, the author's family raises money for JDRF through the ``JDRF One Walk'' fundraiser.  The name of our fundraising team is ``The Fighting Pancreases of Zachary University'', and John and Xavier Golden designed the logo (see Figure 1).  Also in Figure 1 are the front and back views of a 3D-printed logo pendant, and the method to create the pendant is the topic of this paper.

\begin{figure}[h]
\centering
\includegraphics[width=0.29 \textwidth]{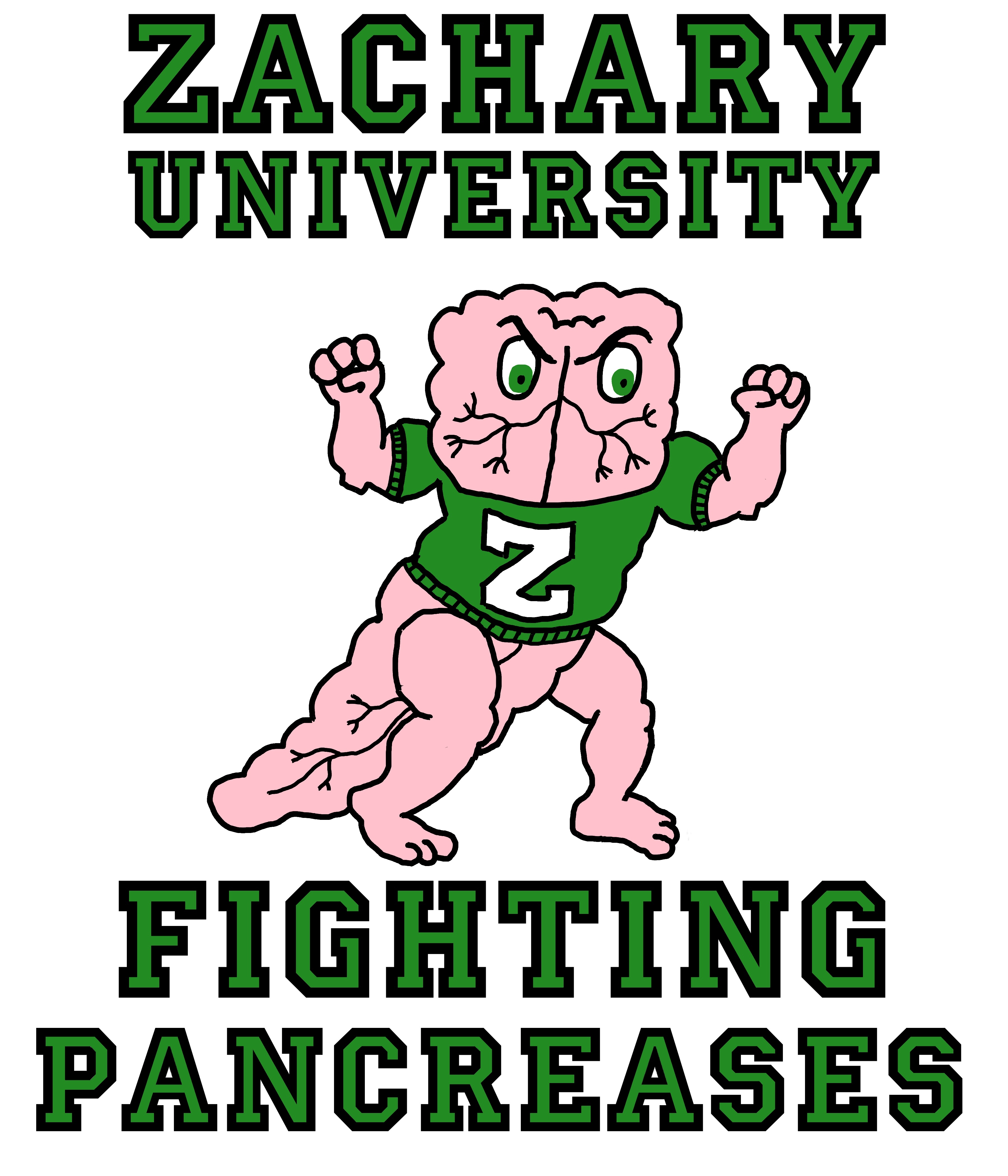}
\includegraphics[width=0.34 \textwidth]{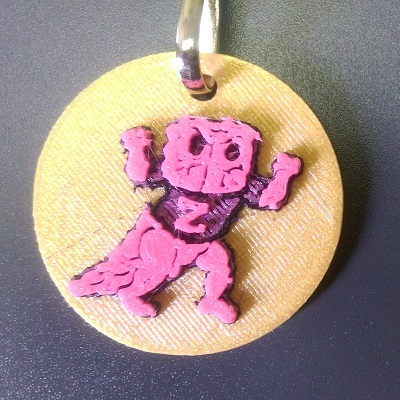}
\includegraphics[width=0.34 \textwidth]{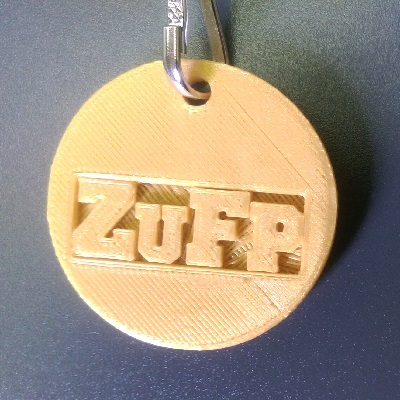}
\caption*{Figure 1: The Fighting Pancreas character, designed by John and Xavier Golden -- original image and 3D-printed pendant.}
\end{figure}

\section{Overview of the Design of the Pendant}

The pendant is constructed using Mathematica by combining three parts -- a base layer that includes the ``ZUFP'' letters, a simple middle layer, and a top layer that is a three-dimensional rendering of the Fighting Pancreas.  Each layer is bounded by the circle $x^2 + (y + 10)^2 = 4900$ and includes the hole for a clip that is defined by $x^2 + (y + 64)^2 \le 49$.  (When looking at the pendants, the positive $y$ axis is down. This is related to how Mathematica stores the data from the image in a matrix.) For both the base and top layers, pre-processed JPEG files are loaded into Mathematica and the received data are transformed using methods particular to this project.  After the three layers are joined together using Mathematica's \texttt{Show} command, the result is exported an STL file that can be used by a 3D printer.

The code below will create a very large pendant, and the author has found that a good first step in the printer software is to scale the model so that the diameter is 50mm.  To create earrings, a bag tag, or a charm, a smaller diameter can be used.  Also, in the actual printing of the pendant, three different colored filaments are used, and as the author has a single-filament printer, the printing process needs to be paused at appropriate times so that the filament can be changed.  This is straightforward to do on a Makerbot Replicator 2 printer.  The author has found that the Hatchbox filament available on Amazon has worked very well to create the pendants.

Details of the steps to create the pendant are below.

\section{The Base Layer -- Mathematica}

The base layer is defined for $0 \le z \le 6$, and the creation of this part of the pendant begins by loading a JPEG file of the flipped version of the ``ZUFP'' letters (see Figure 2) into Mathematica, using code created in 2013 by Melissa Sherman-Bennett and Sylvanna Krawczyk, who were my research students at our REU program.  (See the URL below for their primer on 3D printing, and other project write-ups.)  Here is the code:

\begin{verbatim}
ClearAll["Global`*"]
DataBase = Import["C:\\data\\3d\\ZUFP\\ZUFP-letters-flip.jpg"];
SizeBase = Import["C:\\data\\3d\\ZUFP\\ZUFP-letters-flip.jpg", "ImageSize"];
GrayBase = ColorConvert[Image[DataBase, "Real"], "Grayscale"];
\end{verbatim}

These lines of code load a JPEG into a matrix \texttt{DataBase} of grayscale values (on a 0 to 1 scale).  Specifically, the first line clears the memory in Mathematica, and the next lines load and process the image file from the directory \texttt{C:$\backslash$data$\backslash$3d$\backslash$ZUFP}.  (The author has found that keeping all files in a directory on the C: Drive on a PC is wise, because the eventual STL file will be at least 20MB, and trying to access a network drive with this size file slows processing considerably.)  We use the \texttt{ColorConvert} command because the Import command creates three matrices (for colors red, blue, and green), even if the JPEG file is greyscale.

\begin{figure}[h]
\centering
\includegraphics[width=0.4 \textwidth]{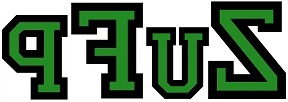}
\includegraphics[width=0.6 \textwidth]{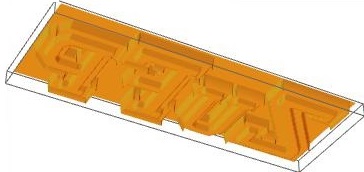}
\caption*{Figure 2: 2D image with flipped letters \& 3D rendering of the letters.}
\end{figure}

Next, we define a function that we are going to apply to every entry in \texttt{GrayBase}.  The basic idea is that the interior of the four letters will be printed from the build plate up, and the other regions will not print at all, until we get to $z = 5$. In order to define this function, it is important to first know what the \emph{luminance} is for each part of the ``ZUFP'' JPEG.  The \emph{luminance} of a color is a measure of the light intensity of the color.  If a color is converted to grayscale, either on a 0 to 1 scale or a 0 to 255 scale, the luminance is the same as the grayscale number.  In \emph{Microsoft Paint}, the luminance of a color in an image can be found by first using the ``Color Picker'' tool to define "Color 1'', and then using ``Edit colors'' to see what the luminance is.  (Note that the Red-Green-Blue numbers for the color can also be found this way, too.)  For the green parts of the letters, the luminance is 82 out of 255 (about 0.32 on the 0-1 scale), and the RGB numbers are 37-140-35 (each out of 255).

To accomplish this, we apply a function called \texttt{BoundBase}.  We are going to transform the data in this way:  the green areas will print to a thickness of 5, the black areas to a thickness of 3, and the white areas will not print at all.  The way to read the following function is first there is a test that we are in the black area of the letters ($x < 0.1$, very dark).  If that test fails, then we go to another test to see if we are in the green area of the letters ($x < 0.7$, which will definitely include 0.32).  If that test also fails, we must be in the white section of the JPEG.  Here is the function:

\begin{verbatim}
BoundBase[x_] := If[x < 0.1, 3, If[x < 0.7, 5, 0]];
\end{verbatim}

The next command in Mathematica is to apply this function to the image and convert to a matrix (or table).  The defined variation of $i$ and $j$ in this command takes into account that the image is stored in memory backwards from what we want to print.  (The \texttt{ArrayPad} command is used to include some extra 0's around the boundary of the matrix, which is sometimes useful to define a boundary.  It might not be necessary for this particular design.)

\begin{verbatim}
ArrayBase = ArrayPad[Table[BoundBase[ImageData[GrayBase][[i, j]]],
     {i, 1, SizeBase[[2]], 1}, {j,  SizeBase[[1]], 1, -1}], {1, 1}, 0];
\end{verbatim}

Now, we can generate an image in Mathematica of the letters, using the important \texttt{ListPlot3D}.  This command is very useful to convert an array of (two-variable) function values into a surface, and it tends to smooth the surface in a way that turns out well using a 3D-printer.  First, here is the command:

\begin{verbatim}
MyBase = ListPlot3D[ArrayBase, DataRange -> {{-56, 56}, {-20, 20}},
    BoxRatios -> Automatic,
    RegionFunction -> Function[{x, y, z}, x^2 + (y + 10)^2 < 4900],
    Mesh -> None, PlotStyle -> {RGBColor[0.85, 0.62, 0.125]},
    BoundaryStyle -> None, Axes -> False];
\end{verbatim}

The \texttt{DataRange} option defines the rectangle that will be the domain of our function of two variables.  This rectangle has to be inside the boundary circle described above, and placed in the appropriate position on the disk.  The \texttt{BoxRatios} option will keep the scaling on the three variable the same, while the \texttt{RegionFunction} option keeps us defined within our disk.  The \texttt{Mesh} and \texttt{BoundaryStyle} command will eliminate any gridlines and edgelines, which is nice for creating images in Mathematica for manuscripts and talks (see Figure 6 at the end of this paper), and the \texttt{PlotStyle} command will color this part of the image in Mathematica a nice gold color (but has no bearing on the choice of filament that you use).  The \texttt{Axes} option can be toggled to true for troubleshooting.  Asking Mathematica to show \texttt{MyBase} gives us the image in Figure 2.

We now need to surround this rectangle with a part of a cylinder in order to finish the first layer of the pendant.  We can define this cylinder using another important Mathematica command for 3Dprinting:  \texttt{RegionPlot3D}.  Here are the commands to create the rest of the base layer:

\begin{verbatim}
MyBase1 = RegionPlot3D[x^2 + (y + 10)^2 < 4900 && x^2 + (y + 64)^2 > 49
    && z >= 0 && z <= 6 && Max[Abs[x], 2.8*Abs[y]] > 55,
    {x, -80, 80}, {y, -80, 80}, {z, -80, 80},
    PlotPoints -> 100, PlotRange -> All, Axes -> True,
    BoxRatios -> Automatic, Mesh -> None,
    PlotStyle -> {RGBColor[0.85, 0.62, 0.125]}, BoundaryStyle -> None];
Show[MyBase1, MyBase]
\end{verbatim}

The key aspects of this code is that \&\& means \emph{and} in Mathematica, the \texttt{Max} section defines the region outside of the letters rectangle (and this is a place where adjustments need to be made depending on the size of the letters JPEG), and the larger the \texttt{PlotPoints} number, the smoother the boundary of the disk, but the longer the code takes to run and the larger the STL file.  The \texttt{Show} command gives us the image in Figure 3.

\begin{figure}[h]
\centering
\includegraphics[width=0.5 \textwidth]{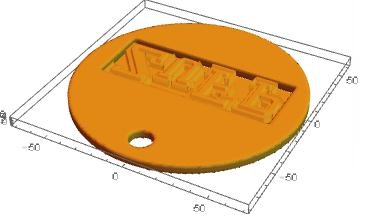}
\caption*{Figure 3: The base layer, with axes on.}
\end{figure}

\section{The Middle Layer -- Mathematica}

The middle layer is simply a cylinder with a hole punched out for the clip. In the code, the middle layer is actually defined to be $5 \le z \le 8$, even though the base layer rises to $z = 6$ and the top layer starts at $z = 7$.  The author has found it useful to have a middle layer and to have the layers intersect -- it seems to help with both the virtual rendering in Mathematica and with the 3D printing -- but an interested reader can experiment with defining the middle layer as $6 \le z \le 7$ or eliminating it altogether.  Here is the code for the middle layer:

\begin{verbatim}
MyMiddle = RegionPlot3D[x^2 + (y + 10)^2 < 4900 && x^2 + (y + 64)^2 > 49
    && z >= 5 && z <= 8, {x, -80, 80}, {y, -80, 80}, {z, -80, 80},
    PlotPoints -> 100, PlotRange -> All, Axes -> True,
    BoxRatios -> Automatic, Mesh -> None,
    PlotStyle -> {RGBColor[0.85, 0.62, 0.125]}];
\end{verbatim}

We can now start to see what happens when we combine layers:

\begin{verbatim}
Show[MyMiddle, MyBase1, MyBase, Axes -> None, Boxed -> False]
\end{verbatim}

Rotating the output a bit will give us the image in Figure 4.

\begin{figure}[h]
\centering
\includegraphics[width=0.4 \textwidth]{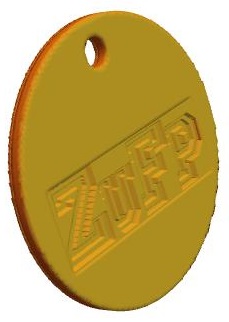}
\caption*{Figure 4: The base and middle layers, with axes and bounding box off.}
\end{figure}

A final note on this part of the design:  in Figure 4, the interiors of the letters are depressed into the pendant, but in the 3D print in Figure 1, the region outside of the letters is depressed.  The interested reader can explore why this is the case.

\section{The Top Layer -- Mathematica}

The creation of the top layer uses many of the idea described above.  However, there are a couple key points to make.  First, if an image, particularly a logo, is too complicated, it will not print well, so it may be necessary to pre-process the image in some way so that there are larger areas of the same color or grayscale.  Second, because the code is going to transform ranges of grayscales to constant heights, it may be necessary or wise to change the colors in certain sections of an image.  (No one needs to know how you have adjusted the JPEG before it is read into Mathematica.)  For this particular pendant, in Figure 5 is the JPEG that the author actually used to create the pendant.  Note that the ``Z'' has been re-colored to match the skin of the Fighting Pancreas.  The author also experimented with more significant recoloring, also in Figure 5.

\begin{figure}[h]
\centering
\includegraphics[width=0.4 \textwidth]{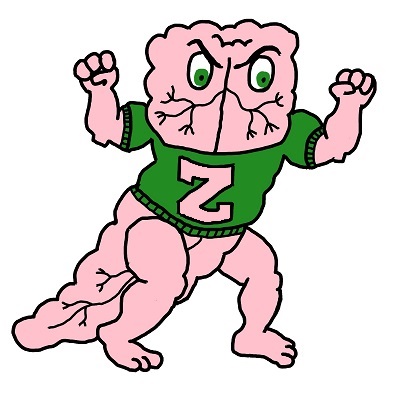}
\includegraphics[width=0.4 \textwidth]{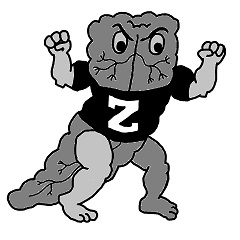}
\caption*{Figure 5: Left -- JPEG for the pendant; Right -- experimental JPEG.}
\end{figure}

The following code creates the top layer of the pendant:

\begin{verbatim}
DataTop = Import["C:\\data\\3d\\ZUFP\\ZUFP-logo-color-small2.jpg"];
SizeTop = Import["C:\\data\\3d\\ZUFP\\ZUFP-logo-color-small2.jpg", "ImageSize"];
GrayTop = ColorConvert[Image[DataTop, "Real"], "Grayscale"];
BoundTop[x_] := If[x > 0.9, 7, If[x < 0.1, 11, If[x < 0.6, 11, 15]]];
ArrayTop = ArrayPad[Table[BoundBase[ImageData[GrayTop][[i, j]]],
    {i, 1, SizeTop[[2]], 1}, {j,  SizeTop[[1]], 1, -1}], {1, 1}, 0];
MyColor[z_] := Piecewise[{{Black, 4 < z <= 10},
    {RGBColor[0/255, 102/255, 0/255 ], 10 < z <= 13}, {Pink, z > 13}}];
MyTop = ListPlot3D[ArrayTop, DataRange -> {{-55, 50}, {-55, 50}},
    Axes -> True, BoxRatios -> Automatic,
    RegionFunction -> Function[{x, y, z}, x^2 + (y + 10)^2 < 4900], Mesh -> None,
    ColorFunction -> MyColor, ColorFunctionScaling -> False, BoundaryStyle -> None];
\end{verbatim}

Of note in this code is that the \texttt{BoundTop} function defines three heights, one for white in the JPEG, the second for pink, and a third for both green and black.  Because of the way the green and black are near each other in the JPEG, the author found that the 3D print turns out better if they are just the same height.  The author also chose to make the ``Z'' the tallest height, which is why it is re-colored in the JPEG that was used.  Also of note is the \texttt{MyColor} function, which colors the 3D rendering in a nice way in Mathematica.

Finally, we can generate the overall image in Mathematica, and then export that image as an STL file.  The \texttt{Export} command takes a minute or two to run:

\begin{verbatim}
MySTL = Show[MyMiddle, MyBase1, MyBase, MyTop, Axes -> None, Boxed -> False]
Export["C:\\data\\3d\\ZackU\\ZackU-Pendant-large2.stl", MySTL]
\end{verbatim}

We can also create some nice images to see what is going on.  The following code is modified from code found online at the useful Stackexchange site\footnote{\url{http://mathematica.stackexchange.com/questions/3759/autorotating-3d-plots}}:

\begin{verbatim}
SnapShot[j_] =Show[MySTL, ViewVertical -> {0, -1, 0}, ViewVector ->
    {RotationMatrix[2 j Pi/40, {0, 1, 0}].
    ({0, -10, 200} - {0, -10, 0}) + {0, -10, 0}, {0, -10, 0}},
    SphericalRegion -> True, ViewAngle -> 50 Degree];
GraphicsRow[{SnapShot[0], SnapShot[5], SnapShot[10], SnapShot[20]},ImageSize -> Large]
For[i = 1, i < numFrames + 1, i++,
    Export[ToString[
    StringForm["C:\\data\\3d\\ZUFP\\Frames\\ZUFrame``.gif",i + 100]], SnapShot[i]]]
\end{verbatim}

The first command allows to create frames for different views of the pendant.  The second command creates the following:

\begin{figure}[h]
\centering
\includegraphics[width=\textwidth]{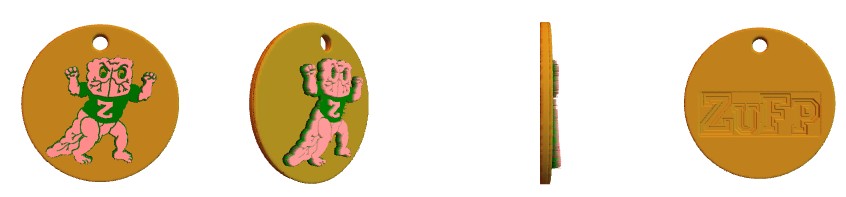}
\caption*{Figure 6: The ZUFP pendant, from various angles.}
\end{figure}

The final command will generate 40 frames of the pendant rotating. These frames can be used to create an animated gif.  Important note:  if the reader plans to run the \texttt{Export} command a second time, delete all the previously-generated GIF ``ZUFrame'' files first, or Mathematica will crash.

\section{Printing the Pendant}

Now that the STL file is created, it can be printed.  As noted above, to create a reasonably-sized pendant, the model should be scaled to something like 50mm diameter, and other scales can be used to create, earrings, etc.  The use of multiple colored filaments is an added challenge, and a test print would be necessary to determine when to pause the printer and switch the filament.  On a Makerbot Replicator 2, using no Raft, 15\% infill, 0.2mm resolution, the author found the following worked:  start with gold filament until the print is 71\% complete.  Then, switch to black or green to start on the Fighting Pancreas, continuing to 86\% for the print in Figure 1.  Note that for that print, there is one layer of the shirt that is pink, giving an interesting mix of pink and black.  Using 87\%, there is no last-minute shading of the shirt.  Complete the arms and legs in pink.

These percentages are developed from trial-and-error.  If the percentages were computed by the Makerbot based on $z$ height, because of the way the pendant is designed, we could determine the percentages \emph{a priori}.  However, the percentages are based on either time remaining or volume remaining (which should be close to each other), so while it would be an interesting calculus problem to find the percentages \emph{a priori}, it may just be easier to do a test print and observe carefully where to change the filament.  Printing one pendant takes about 25 minutes on the Makerbot Replicator 2.

In Figure 7 are other pendants that have been created using the procedure described in this document, as well as some other pre-processing and Mathematica tricks that may be shared in a future paper or talk.  STL files for two of these pendants, and others not indicated here, can be found at \url{http://www.thingiverse.com/aboufade/designs}.

\begin{figure}[h]
\centering
\includegraphics[width=0.26 \textwidth]{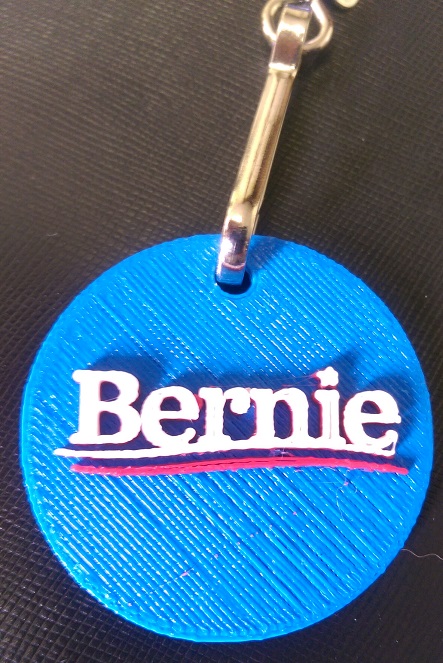}
\includegraphics[width=0.30 \textwidth]{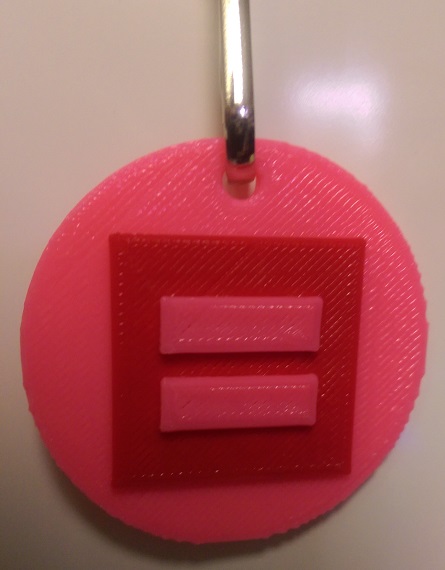}
\includegraphics[width=0.29 \textwidth]{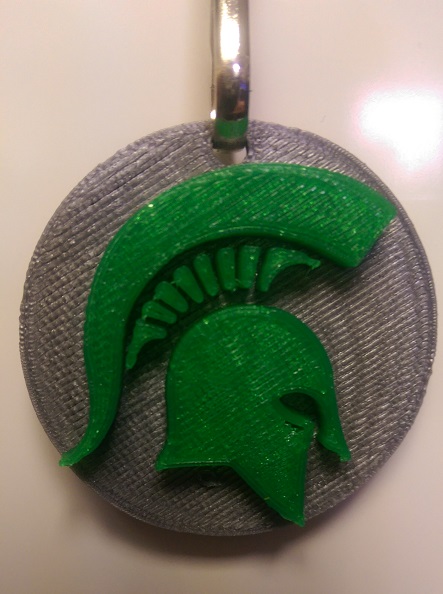}
\caption*{Figure 7: Other pendants by the author.}
\end{figure}

\section{Appendix: URLs}
\begin{itemize}
\item JDRF home page \url{http://www.jdrf.org/}
\item Prof.~Aboufadel's 3D printing page \url{http://sites.google.com/site/aboufadelreu/Profile/3d-printing}
\end{itemize}

\end{document}